\documentclass[a4paper]{article}

\usepackage[UKenglish]{babel}
\usepackage[style=alphabetic,maxbibnames=99]{biblatex}
\usepackage{amsmath}
\usepackage{mathtools}
\usepackage{amsthm}
\usepackage{thmtools}
\usepackage{csquotes}
\usepackage[inline]{enumitem}
\usepackage{nicefrac}
\usepackage{fonttable}
\usepackage{tikz}
\usepackage[hidelinks]{hyperref}
\usepackage[capitalize,noabbrev]{cleveref}
\usepackage{newunicodechar}
\usepackage[final]{microtype}
\usepackage{authblk}

\usepackage[charter]{mathdesign}
\usepackage[T1]{fontenc}

\usepackage[draft]{fixme}

\fxusetheme{color}

\newtheorem{theorem}{Theorem}
\newtheorem{lemma}[theorem]{Lemma}
\newtheorem{proposition}[theorem]{Proposition}
\newtheorem{corollary}[theorem]{Corollary}
\newtheorem{question}[theorem]{Question}
\newtheorem*{mainquestion}{Main question}
\theoremstyle{definition}
\newtheorem{definition}[theorem]{Definition}
\newtheorem{construction}[theorem]{Construction}

\theoremstyle{remark}
\newtheorem{remark}[theorem]{Remark}

\usetikzlibrary{
	cd,
	positioning,
	arrows,
 	arrows.meta,
	calc,
	external,
	math,
	decorations.pathmorphing,
	decorations.markings,
	calligraphy
}
\tikzset{
	point/.style={draw=black,fill=black,opacity=1,circle,outer sep=0pt,inner sep=0,minimum size=2},
	dot/.style={draw=black,fill=black,opacity=0,circle,outer sep=0pt,inner sep=0}
}
\tikzset{
	negated/.style={
        decoration={markings,
            mark= at position 0.5 with {
                \node[transform shape] (tempnode) {$\backslash$};
            }
        },
        postaction={decorate}
    }
}
\tikzset{
	brace/.style={decoration={calligraphic brace,amplitude=5pt}, decorate, line width=1.25pt}
}

\newunicodechar{Λ}{\ensuremath{\Lam}}

\newcommand{\Pad}{\mathrm{Pad}}
\newcommand{\suc}{\mathsf{succ}}
\newcommand{\In}{\mathrm{In}}
\DeclareMathOperator{\Road}{\mathrm{Road}}
\DeclareMathOperator{\Branch}{\mathrm{Branch}}

\newcommand{\R}{\ensuremath{\mathbb{R}}}
\newcommand{\Q}{\ensuremath{\mathbb{Q}}}
\newcommand{\Z}{\ensuremath{\mathbb{Z}}}

\DeclareMathOperator{\dom}{\mathrm{dom}}
\DeclareMathOperator{\ran}{\mathrm{ran}}

\DeclareMathOperator{\height}{\mathrm{height}}
\newcommand{\restr}{\mathord{\upharpoonright}}

\newcommand{\cont}{\mathfrak{c}}
\DeclareMathOperator{\cf}{\mathrm{cf}}
\newcommand{\ZFC}{\ensuremath{\mathrm{ZFC}}}
\newcommand{\CH}{\ensuremath{\mathrm{CH}}}
\newcommand{\SH}{\ensuremath{\mathrm{SH}}}
\newcommand{\MA}{\ensuremath{\mathrm{MA}}}

\DeclareMathOperator{\dS}{\downarrow}
\DeclareMathOperator{\uS}{\uparrow}
\makeatletter
\newcommand{\@usstar}[1]{{\uparrow}\left(#1\right)}
\newcommand{\@usnostar}[1]{{\uparrow}(#1)}
\newcommand{\us}{\@ifstar{\@usstar}{\@usnostar}}
\newcommand{\@dsstar}[1]{{\downarrow}\left(#1\right)}
\newcommand{\@dsnostar}[1]{{\downarrow}(#1)}
\newcommand{\ds}{\@ifstar{\@dsstar}{\@dsnostar}}
\makeatother

\newcommand{\lam}{\lambda}
\newcommand{\Lam}{\ensuremath\Lambda}

\newcommand{\Ra}{\Rightarrow}
\newcommand{\Lra}{\Leftrightarrow}
\newcommand{\es}{\varnothing}
\newcommand{\sse}{\subseteq}
\newcommand{\partto}{\rightharpoonup}

\newcommand{\wh}{\widehat}
\newcommand{\bb}{\mathbb}
\newcommand{\mc}{\mathcal}
\newcommand{\mr}{\mathrm}
\newcommand{\defeq}{\vcentcolon=}
\DeclarePairedDelimiter{\ab}{\langle}{\rangle}
\DeclarePairedDelimiter{\abs}{|}{|}

\renewcommand{\leq}{\leqslant}
\renewcommand{\geq}{\geqslant}

\newcommand{\contradiction}{\noindent
	\begin{tikzpicture}[x=0.4ex,y=0.4ex]
		\draw[line width=.15ex] (0,0) -- (1,2) -- (0,2) -- (1,4)
		(0.95,0.32) -- (0,0) -- (-0.32,0.95);
	\end{tikzpicture}\hspace*{0.2em}
}

\title{On the continuous gradability of the cut-point orders of \texorpdfstring{\R}{R}-trees}
\author[*]{Sam Adam-Day}
\affil[*]{Mathematical Institute, University of Oxford, Andrew Wiles Building, Radcliffe Observatory Quarter, Woodstock Road, Oxford, OX2 6GG, United Kingdom; \href{mailto:adamday@maths.ox.ac.uk}{\nolinkurl{adamday@maths.ox.ac.uk}}}
\date{}

\begin{document}

	\maketitle

	\abstract{
		An \R-tree is a certain kind of metric space tree in which every point can be branching. Favre and Jonsson posed the following problem in 2004: can the class of orders underlying \R-trees be characterised by the fact that every branch is order-isomorphic to a real interval? In the first part, I answer this question in the negative: there is a `branchwise-real tree order' which is not `continuously gradable'. In the second part, I show that a branchwise-real tree order is continuously gradable if and only if every well-stratified subtree is \R-gradable. This link with set theory is put to work in the third part answering refinements of the main question, yielding several independence results. For example, when $\kappa \geq \cont$, there is a branchwise-real tree order which is not continuously gradable, and which satisfies a property corresponding to $\kappa$-separability. Conversely, under Martin's Axiom at $\kappa$ such a tree does not exist.
	}

	\renewcommand{\thefootnote}{}
	\renewcommand{\thefootnote}{\arabic{footnote}}
	\addtocounter{footnote}{-2}

	\listoffixmes


\section{Introduction}
\label{sec:intro}

An \R-tree is to the real numbers what a graph-theoretic tree is to the integers. Formally, let $\ab{X,d}$ be a metric space. An \emph{arc} between $x,y \in X$ is the image of a topological embedding $r \colon [a,b] \to X$ of a real interval such that $r(a) = x$ and $r(b) = y$ (allowing for the possibility that $a=b$). The arc is a \emph{geodesic segment} if $r$ can be taken to be an isometry. The metric space $\ab{X,d}$ is an \emph{\R-tree} if between any two points $x,y \in X$ there is a unique arc, denoted $[x,y]$, which is also a geodesic segment.

Note that any tree in the graph-theoretic sense can be viewed as an \R-tree via its so-called `geometric realisation'. Indeed, let $V$ be a (possibly infinite) set of vertices and let $E$ be a symmetric binary relation on $V$, such that $G = \ab{V,E}$ is a connected graph without cycles. Then $G$ can be realised as an \R-tree by taking $V$ as a discrete set of points and adding a copy of the unit interval $(0,1)$ between $u,v \in V$ whenever $\ab{u,v} \in E$.

\begin{figure}
	\begin{equation*}
		\begin{tikzpicture}
			\begin{scope}[yscale=-1]
				\node at (-1,0) {(a)};
				\draw[gray] (0,0) -- (6.5,0);
				\draw[very thick, olive, every node/.style={point,olive}]
					(1,1.5) node {} -- (1,0.25) node {};
				\draw[very thick, blue, every node/.style={point,blue}]
					(2,-1.1) node {} -- (2,0) -- (3,0) -- (3,-0.75) node {};
				\draw[very thick, red, every node/.style={point,red}]
					(3.9,-0.4) node {} -- (3.9,0) -- (4.4,0) -- (4.4,0.9) node (pa) {};
				\draw[very thick, orange, every node/.style={point,orange}]
					(5.5,0.3) node {} -- (5.5,-1) node {};
				\node[right=1pt of pa] {$p$};
			\end{scope}
			\begin{scope}[yshift=-160]
				\node at (-1,1) {(b)};
				\draw[gray] (4,0) edge +(135:4.4) edge +(45:2.1);
				\draw[very thick, olive, every node/.style={point,olive}]
					($(4,0) + (135:3.4) + (165:0.35)$) node {} -- ++(165:1.25) node {};
				\draw[very thick, blue, every node/.style={point,blue}]
					($(4,0) + (135:1.4) + (105:0.75)$) node {} -- ++(105:-0.75) -- ++(135:1) -- ++(105:1.1) node {};
				\draw[very thick, red, every node/.style={point,red}] 
					(4,-0.9) node (pb) {} -- (4,0) -- ++(135:0.5) -- ++(105:0.4) node {};
				\draw[very thick, orange, every node/.style={point,orange}]
					($(4,0) + (45:1.1) + (15:0.3)$) node {} -- ++(15:-0.3) -- ++(75:1) node {};
				\node[right=1pt of pb] {$p$};
			\end{scope}
		\end{tikzpicture}
	\end{equation*}
	\caption{(a) Some example geodesic segments in the train track \R-tree. (b) The cut-point order on this tree with root $p$, showing the images of the example geodesic segments.}
	\label{fig:example geodesics on non-simplicial tree}
\end{figure}
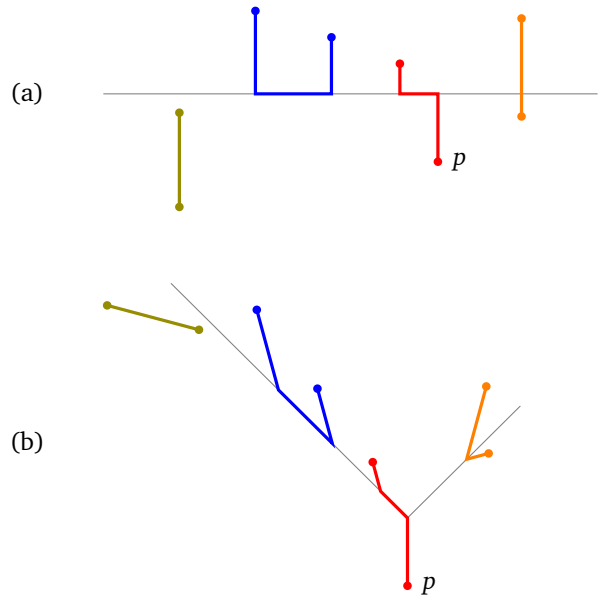

The class of \R-trees however is much more general than this. Consider the following example of an \R-tree which does not arise in this fashion. Let $X$ be the space obtained by taking the real plane $\R^2$, and designating each point on the $x$-axis as a `train station' and each vertical line $\{x\} \times \R$, as well as the $x$-axis, as a `train track'. We define a new metric on $X$: to travel between two points in the real plane, one must travel along the tracks, potentially passing through train stations. In the resulting \R-tree, the removal of any point on the $x$-axis leaves exactly $4$ connected components. See \cref{fig:example geodesics on non-simplicial tree}(a) for some example geodesic segments in this tree.

Now, \R-trees play an important role in geometric group theory, and are interesting objects in their own right. (See \cite{bestvina1997real,mayer92} and the references contained in \cite{Fabel15}.) In this paper, we are interested in the underlying order structure of \R-trees. Given an \R-tree $X$ together with a designated point $p \in X$, the \emph{cut-point order} on $X$ with root $p$ is defined by, for $x, y \in X$:
\begin{equation*}
	x \leq y \quad\Lra\quad [p,x] \sse [p,y]
\end{equation*}
See \cref{fig:example geodesics on non-simplicial tree}(b) for a picture of a cut-point order on the train track \R-tree example. 

In \S 3.1 of \cite{favre2002valuative}, Charles Favre and Mattias Jonsson investigate this order structure, one of their aims being an order-theoretic characterisation of structures arising in this way. I replicate their definitions here, with a few changes of terminology.

\begin{definition}\label{def:tree order}
	A \emph{tree order} is a partial order $X$ such that the following conditions hold.
	\begin{enumerate}[label=(TO\arabic*), leftmargin=*, labelindent=3pt]
		\item\label{item:lo; def:tree order}
			For every $x \in X$ the set $\ds x = \{y \in X \mid y \leq x\}$ is a linear order.
		\item\label{item:root; def:tree order}
			$X$ has a minimum element, its \emph{root}.
	\end{enumerate}
	A \emph{branch} in $X$ is a maximal linearly ordered subset.
\end{definition}

Partial orders satisfying \ref{item:lo; def:tree order} are sometimes called `psuedotrees' (see for instance \cite{nikiel1989topologies}).

\begin{definition}\label{def:brot}
	A \emph{branchwise-real tree order} is a tree order $X$ subject to the following extra conditions.
	\begin{enumerate}[label=(BR\arabic*), leftmargin=*, labelindent=3pt]
		\item\label{item:interval; def:brot}
			Every branch is order-isomorphic to a real interval.
		\item\label{item:meet-semilattice; def:brot}
			$X$ is a meet-semilattice; that is, any two points $x, y \in X$ have a greatest lower bound $x \wedge y$, their \emph{meet}.
	\end{enumerate}
\end{definition}

Favre and Jonsson call such objects `non-metric trees'.

\begin{remark}
	In fact, Favre and Jonsson's definition is equivalent to \ref{item:lo; def:tree order} + \ref{item:root; def:tree order} + \ref{item:interval; def:brot}. They erroneously claim that \ref{item:meet-semilattice; def:brot} follows from the rest using the completeness of the real line \cite[p.~45]{favre2002valuative}. This is incorrect in light of the following counterexample. Let $X$ be the interval $[0,1)$, together with two incomparable copies of the element $1$ sitting on top, as in the following diagram.
	
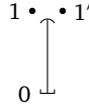
\begin{equation*}
		\begin{tikzpicture}
			\draw[[-)] (0,0) node [label=left:$0$] {} -- (0,1);
			\node[point] at (-0.2,1.1) [label=left:$1$] {};
			\node[point] at (0.2,1.1) [label=right:$1'$] {};
		\end{tikzpicture}
	\end{equation*}
	Then $1$ and $1'$ have no common meet, and so $X$ is not a branchwise-real tree order.
\end{remark}

\begin{definition}
	Let $P$ and $Q$ be partial orders. A \emph{$Q$-grading} of $P$ is a strictly monotonic map $f \colon P \to Q$. That is, whenever $x < y$ in $P$ we have $f(x) < f(y)$.
\end{definition}

\begin{definition}\label{def:continuous grading}
	Let $X$ be a branchwise-real tree order. An \R-grading $\ell \colon X \to \R$ is \emph{continuous} if for any $x < y$ in $X$, letting $[x,y] \defeq \{z \in X \mid x \leq z \leq y\}$, the restriction: 
	\begin{equation*}
		\ell\restr{[x,y]} \colon [x,y] \to [\ell(x),\ell(y)]
	\end{equation*}
	is an order-isomorphism. I will usually drop the `\R' and call such functions \emph{continuous gradings}. Say that $X$ is \emph{continuously gradable} if it admits a continuous grading. 
\end{definition}

In \cite{favre2002valuative} these functions are called `parametrizations'. 

We can now state the order-theoretic characterisation of \R-tree cut-point orders which Favre and Jonsson obtained.

\begin{theorem}\label{res:cut-point para brots}
	The class of \R-tree cut-point orders is exactly the class of continuously gradable branchwise-real tree orders.
\end{theorem}

\begin{figure}[t]
	\begin{equation*}
		\begin{tikzpicture}
			\begin{scope}[every node/.style={point}]
				\node (p) at (0,0) [label=below:$p$] {};
				\node (m) at (0,2) {};
				\node (x) at (-2,4) [label=above left:$x$] {};
				\node (y) at (2,4) [label=above right:$y$] {};
			\end{scope}
			\node (ml) at ($(m) + (1.5,-0.3)$) {$x \wedge y$};
			\draw[->] (ml)  -- (m);
			\draw (p) -- (m) -- (x);
			\draw (m) -- (y);
			\begin{scope}[every path/.style={{Bar-Bar}}]
				\draw[green!50!black] ($(p)+(-0.6,0)$) -- node[left] {$\ell(x \wedge y)$} ($(m)+(-0.6,0)$);
				\draw[blue] ($(p)+(-0.2,0)$) -- ($(m)+(-0.185,-0.077)$) -- node[below left] {$\ell(x)$} ($(x)+(-0.141,-0.141)$);
				\draw[blue] ($(p)+(0.2,0)$) -- ($(m)+(0.185,-0.077)$) -- node[below right] {$\ell(y)$} ($(y)+(0.141,-0.141)$);
				\draw[red] ($(x)+(0.141,0.141)$) -- ($(m) + (0,0.283)$) node[above=0.5] {$d(x,y)$} -- ($(y)+(-0.141,0.141)$);
			\end{scope}
		\end{tikzpicture}
	\end{equation*}
	\caption{The railroad track equation}
	\label{fig:railroad track equation}
\end{figure}
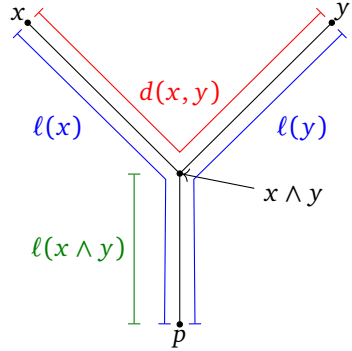

\begin{proof}
	See \cite[p.~50]{favre2002valuative}. Given an \R-tree $X$ and $p \in X$, it is straightforward to verify that the cut-point order on $X$ with root $p$ is a branchwise-real tree order. Moreover, we can use the metric $d$ on $X$ to define a continuous grading by:
	\begin{equation*}
		\ell(x) \defeq d(p,x)
	\end{equation*}
	Conversely, given a branchwise-real tree order $X$ and a continuous grading $\ell \colon X \to \R$, we can use $\ell$ to define a metric on $X$ via the `railroad track equation' (see also \cite{mayer92}):
	\begin{equation*}
		d(x,y) \defeq  \ell(x) + \ell(y) - 2\ell(x \wedge y)
	\end{equation*}
	See \cref{fig:railroad track equation} for an illustration of railroad track equation.
\end{proof}

The construction of the metric on the branchwise-real tree order $X$ relies on the continuous grading. It is natural however to imagine that one could always find a continuous grading on a branchwise-real tree order using its real-line-like structure. But Favre and Jonsson state: \textquote[{\cite[p.~47]{favre2002valuative}}]{We do not know if there exists a non-parameterizable nonmetric tree}. This challenge is taken up in the present paper, in which I answer the following main question.

\begin{mainquestion}[\citeauthor{favre2002valuative}]\label{q:main}
	Is every branchwise-real tree order continuously gradable?
\end{mainquestion}

At its heart, this question asks about the existence of a local-global connection for these tree orders. Each branch of a branchwise-real tree order $X$ embeds into \R, and the question is whether these `local' embeddings can be combined consistently into one `global' continuous grading $X \to \R$. In \cref{sec:non-para} I will show that the answer to the main question is in fact `No'. In other words, there is a branchwise-real tree order in which these local \R-embeddings do not combine into one complete continuous grading:

\begin{theorem}
	There is a branchwise-real tree order which is not continuously gradable.
\end{theorem}

In order to construct our non-continuously-gradable branchwise-real tree order, we will step into the realm of set-theoretic trees (here referred to as `well-stratified trees', for the sake of clarity). For definitions of the set-theoretic concepts used in remainder of this introduction the reader is referred to \cref{sec:set theory}. Baumgartner, Laver and Gavin showed in 1970 that there exists a well-stratified tree, all of whose branches are countable, which admits no \R-grading (throughout `countable' means `finite or countably infinite'). The construction of this tree is given in \cref{sec:non-para}. Once we have such a tree $T$, we take what is known as its `road space': each successor node of $T$ is replaced by a copy of the real interval $(0,1]$. The resulting tree order is then branchwise-real and has no continuous grading, as required.

The remainder of the paper is concerned with developing refinements of this technique. In \cref{sec:inside}, the link between well-stratified trees and branchwise-real tree orders is strengthened by proving this following result. 

\begin{theorem}
	A branchwise-real tree order $X$ is continuously gradable if and only if every well-stratified subtree $T \sse X$ is \R-gradable.
\end{theorem}

The proof works by constructing an increasing sequence $(T_n)_{n \in \omega}$ of `approximating' well-stratified subtrees of $X$, whose union is dense in the interval topology; that is, for any $x < y$ in $X$ there is $n \in \omega$ and $z \in T_n$ with $x < z < y$. This technique allows for the application of set-theoretic methods to problems involving branchwise-real tree orders.

The third part of the paper focuses on answering certain refinements of the main question, some of which turn out to be independent of \ZFC. Can we obtain a non-continuously-gradable branchwise-real tree order satisfying certain additional properties? Or, conversely, which properties of branchwise-real tree orders entail continuous gradability?

\Cref{sec:branching nodes} considers a property motivated from the study of \R-trees. Those \R-trees which appear in applications are often separable (see \cite{bestvina1997real}). The first result is that an \R-tree is separable if and only if every cut-point order, regardless of root point, contains at most countably many branching nodes and at most countably many maximal terminal segments isomorphic to a non-trivial real interval. Such orders will be called \emph{countably wispy}. In contrast to the original, the corresponding refinement of the main question has a positive answer, as follows.

\begin{theorem}\label{intro-res:count-wispy cont gradable}
	Every countably wispy branchwise-real tree order is continuously gradable.
\end{theorem}

Generalising, we can ask whether a branchwise-real tree order satisfying the natural generalisation of countable wispiness — \emph{$\kappa$-wispiness}, for some uncountable cardinal $\kappa$ — is automatically continuously gradable. Analogously to the countable case, $\kappa$-wispiness corresponds to $\kappa$-separability on \R-trees. In the second part of \cref{sec:branching nodes}, I show that the answer to this question in general is independent of \ZFC\ set theory, as follows. Note that part \ref{item:MA; intro-res:wispiness independence} is a generalisation of \cref{intro-res:count-wispy cont gradable}.

\begin{theorem}\label{intro-res:wispiness independence}\
	\begin{enumerate}[label=(\arabic*)]
		\item\label{item:CH; intro-res:wispiness independence}
			If the Continuum Hypothesis holds then there is an $\aleph_2$-wispy branchwise-real tree order with no continuous grading.
		\item\label{item:MA; intro-res:wispiness independence}
			If Martin's Axiom holds at $\kappa$ then all $\kappa$-wispy branchwise-real tree orders are continuously gradable.
	\end{enumerate}
\end{theorem}

In the final section I consider branchwise-real tree orders satisfying the countable chain condition (ccc). On \R-trees, this condition corresponds to a property which is slightly stronger than separability. The answer to the corresponding refinement of the main question is independent of \ZFC, as follows.

\begin{theorem}
	The Suslin Hypothesis is equivalent to the statement that every ccc branchwise-real tree order has a continuous grading.
\end{theorem}

The countable chain condition generalises to the $\kappa$ chain condition ($\kappa$-cc). The correspondent on \R-trees is a property slightly stronger than $\kappa$-separability. To answer the generalised question `is every $\kappa$-cc branchwise-real tree order continuously gradable?', I define the notion of a `${<}\kappa$-wide $\R$-ungradable tree', which generalises that of a Suslin tree, and show that the existence of such an object is in general independent of \ZFC. Finally, I obtain the following independence result.

\begin{theorem}
	Let $\kappa$ be an uncountable cardinal. There exists a $\kappa$-cc branchwise-real tree order with no continuous grading if and only if there exists a ${<}\kappa$-wide $\R$-ungradable tree.
\end{theorem}

The paper is concluded with a number of open questions.
	\pagebreak

\section{Background definitions}
\label{sec:set theory}

The following sets out the main definitions from combinatorial set theory which will come into play in this paper. For background on these topics, the reader may consult \cite[Ch.~9]{jech} or \cite[\S~III.3, \S~II.5]{kunen}. A \emph{well-stratified tree} $T$ is a tree order in which every $\ds x$ is well-ordered (in a purely set-theoretic context, we would simply say `tree'). For any $x \in X$ its \emph{rank} is the order type of $\ds x \setminus \{x\}$. For any ordinal $\alpha$, the \emph{$\alpha$th level} of $T$, denoted $T(\alpha)$, is the set of elements of rank $\alpha$. The \emph{height} of $T$, denoted by $\height(T)$, is the least $\alpha$ such that $T(\alpha)$ is empty. Say that $T$ is \emph{Hausdorff} if it is a meet-semilattice (i.e.\@ satisfies \ref{item:meet-semilattice; def:brot}).

An \emph{antichain} in a tree order $X$ is a subset $A \sse X$ such that for any distinct $x, y \in A$ we have $x \not\leq y$ and $y \not\leq x$. Say that $X$ has the \emph{countable chain condition} (ccc) if it has no uncountable antichains.

An \emph{Aronszajn tree} is a well-stratified tree of height $\omega_1$, all of whose levels and branches are countable. A \emph{Suslin tree} is a ccc Aronszajn tree. \emph{Suslin's Hypothesis} (\SH) is the statement that there are no Suslin trees. It is well-known that \SH\ is independent of \ZFC\ (see for instance \cite[p.~239--242]{jech}).

A \emph{forcing poset} is a partial order with a maximum element. In this paper, while we think of trees as growing upwards, forcing posets are thought of as extending downwards, as is standard in set theory. For this reason, the notion of a `ccc forcing poset' is defined in the opposite way to the notion of a `ccc tree order'. Let $\bb P$ be any forcing poset. An \emph{antichain} in $\bb P$ is a subset $A \sse \bb P$ such that for any distinct $x, y \in A$ there is no $z \in \bb P$ such that $z \leq x,y$. Then $\bb P$ has the \emph{countable chain condition} (ccc) if it has no uncountable antichains.

A subset $D \sse \bb P$ is \emph{dense} if for any $p \in \bb P$ there is $q \leq p$ such that $q \in D$. A subset $F \sse \bb P$ is a \emph{filter} if the following hold.
\begin{enumerate}[label=(F\arabic*), leftmargin=*, labelindent=3pt]
	\item $F$ is non-empty.
	\item $F$ is upwards-closed: if $p \in F$ and $p \leq q$ then $q \in F$.
	\item For any $p,q \in F$ there is $r \in F$ such that $r \leq p,q$.
\end{enumerate}
For $\kappa < \cont$, \emph{Martin's Axiom} for $\kappa$, denoted $\MA_\kappa$, is the statement that for any ccc forcing poset $\bb P$ and any collection $\cal D$ of $\kappa$-many dense subsets of $\bb P$, there is a filter $G \sse \bb P$ which intersects every element of $\cal D$. \emph{Martin's Axiom} is the statement that $\MA_\kappa$ holds for all $\kappa < \cont$. It is well-known that $\MA + \cont = \aleph_\alpha$\ is independent of $\ZFC$, for any regular and uncountable $\aleph_\alpha$ (see \cite[Theorem~16.13 and Theorem~16.16]{jech}).

Let me now indicate the general mathematical conventions which will be followed in this paper. A partial function between sets $X$ and $Y$ will be denoted using the notation $f \colon X \partto Y$. Let $P$ be any partial order. For any $x \in P$ let:
\begin{equation*}
	\us x \defeq \{y \in P \mid y \geq x\}
\end{equation*}
If $S \sse P$ is any subset, let:
\begin{equation*}
	\dS S \defeq \{y \in P \mid \exists x \in S \colon y \leq x\},\quad \uS S \defeq \{y \in P \mid \exists x \in S \colon y \geq x\}
\end{equation*}
For any $x < y$ in $P$, define:
\begin{equation*}
	[x,y] \defeq \{z \in P \mid x \leq z \leq y\}
\end{equation*}
Define the other intervals $[x,y)$, $(x,y]$ and $(x,y)$ analogously. Throughout, unless otherwise specified, `monotonic function' means `strictly monotonic function'. Tuples will be denoted using angle brackets: $\ab{a,b, \ldots}$.

\section{Answering the main question}
\label{sec:non-para}

In this section, I answer the main question by proving the following result.

\begin{theorem}\label{res:non cont grad brto}
	There is a branchwise-real tree order which is not continuously gradable.
\end{theorem}

Our non-continuously-gradable branchwise-real tree order will be constructed by taking the `road space' of a certain well-stratified tree. This notion of the road space was first introduced by Floyd Burton Jones \cite{REMARKSONTHENORMALMOORESPACEMETRIZATIONPROBLEM}.

\begin{definition}
	Let $T$ be a well-stratified tree. The \emph{road space} of $T$, denoted $\Road(T)$, is the partial order obtained by replacing each node on a successor level with a copy of the real interval $(0,1]$, and every other node with a copy of the element $1$. Formally, we can view $\Road(T)$ as the following suborder of the lexicographic product order $T \times (0,1]$:
	\begin{equation*}
		\Road(T) = \{\ab{x,t} \in T \times (0,1] \mid t=1 \text{ or }x \in T(\alpha+1)\text{ for some }\alpha\}
	\end{equation*}
\end{definition}

Note that $T$ embeds canonically in $\Road(T)$ via $x \mapsto \ab{x,1}$.

\begin{lemma}\label{res:road space of cb tree is brot}
	When $T$ is a Hausdorff well-stratified tree with no uncountable branches its road space is a branchwise-real tree order.
\end{lemma}

To prove this, we make use of a basic set-theoretic result, which will reappear often enough to warrant a number.

\begin{lemma}\label{res:countable ordinal in Q}
 	Every countable ordinal $\alpha$ embeds into \Q\ in such a way that limits are preserved.
\end{lemma}

\begin{proof}
	A slick proof makes use of the `forth' part of the classical `back-and-forth method'. Enumerate $\alpha = \{\beta_n \mid n \in \omega\}$. We then build up an embedding $f \colon \alpha \to \Q$ inductively on the enumeration. Once we have $f \restr \{\beta_0, \ldots, \beta_{n-1}\}$, since \Q\ is a dense linear order without endpoints, we can find $f(\beta_{n})$ whose relative position with respect to $f(\beta_0), \ldots, f(\beta_{n-1})$ is the same that of $\beta_n$ with respect to $\beta_0, \ldots, \beta_{n-1}$. To ensure that limits are also preserved in the resulting embedding, it suffices to require in addition that the distance between $f(\beta_n)$ and its immediate successor in $f(\beta_0), \ldots, f(\beta_{n-1})$ with respect to the order on \Q\ (if this successor exists) is less than $\nicefrac 1 {n+1}$.
\end{proof}

\begin{proof}[Proof \cref{res:road space of cb tree is brot}]
	Since $T$ is a meet semilattice, so is $\Road(T)$, and hence \ref{item:meet-semilattice; def:brot} is satisfied. As for \ref{item:interval; def:brot}, take any branch $B$ in $\Road(T)$ with the aim of showing that it is isomorphic to a real interval. Let $B_T$ be the result of restricting $B$ to the canonical embedded copy of $T$ in $\Road(T)$. Then $B_T$ is a branch in $T$, and hence by assumption it is isomorphic to a countable ordinal $\alpha$. Hence by \cref{res:countable ordinal in Q} there is a limit-preserving embedding $\alpha \to \Q$, which then extends to an embedding $B \to \R$, whose image is a real interval.
\end{proof}

Note that if $\ell \colon \Road(T) \to \R$ is a continuous grading, then $\ell$ restricts to an \R-grading of $T$ (via the canonical embedding). Thus, with this lemma, to find a branchwise-real tree order $X$ with no continuous grading it suffices to find a Hausdorff, \R-ungradable well-stratified tree with no uncountable branches. James E.\@ Baumgartner first constructed a well-stratified tree with these properties, using results from Richard Laver and F.\@ Gavin (see \cite{baumthesis}). To understand the construction, let us first examine what it means for a well-stratified tree to be \R-gradable, and relate this to \Q-gradability. First, \Q-gradability is equivalent to the well-known notion of `specialness'.

\begin{definition}
	A well-stratified tree is \emph{special} if it is the union of countably many antichains.
\end{definition}

\begin{lemma}\label{res:Q-grading iff special}
	A well-stratified tree $T$ is \Q-gradable if and only if it is special.
\end{lemma}

\begin{proof}
	See \cite[Lemma~III.5.17]{kunen}. We only need the forwards direction, which works as follows. Let $f \colon T \to \Q$ be a \Q-grading. Let $(q_n)_{n \in \omega}$ be an enumeration of \Q. Then for each, $n \in \omega$, let $A_n \defeq f^{-1}\{q_n\}$. Each $A_n$ is an antichain and $T = \bigcup_{n \in \omega} A_n$.
\end{proof}

It turns out that \R-gradability is equivalent to the \Q-gradability of the nodes of the tree with successor rank.

\begin{definition}
	Let $T(\suc)$ be the suborder of $T$ consisting of the nodes on its successor levels.
\end{definition}

\begin{lemma}\label{res:R-grading iff succ-special}
	Let $T$ be a well-stratified tree. Then $T$ is \R-gradable if and only if $T(\suc)$ is \Q-gradable.
\end{lemma}

This result is to be found in \cite[Ch.~4, Theorem 1(b)]{baumthesis}. According to Baumgartner, it is due to Gavin (unpublished).

\begin{proof}
	Assume that $f \colon T \to \R$ is monotonic. For $x \in T(\suc)$ with immediate predecessor $y$, choose $g(x) \in (f(y),f(x)] \cap \Q$. Then $g \colon T(\suc) \to \Q$ is monotonic. Conversely, assume that $T(\suc)$ is \Q-gradable. By \cref{res:Q-grading iff special} then we have that $T(\suc) = \bigcup_{n \in \omega}A_n$, where each $A_n$ is an antichain. Define:
	\begin{align*}
		f \colon T &\to \R \\
		x &\mapsto \sum \left\{\frac{1}{n^2} \;\middle|\; \exists y \leq x \colon y \in A_n\right\}
	\end{align*}
	Then $f$ is monotonic.
\end{proof}

Given any well-stratified tree $T$, we can obtain it as $T'(\suc)$ of some other tree $T'$ as follows.

\begin{definition}
	Let $T$ be a well-stratified tree. Let $\Pad(T)$ be the result of adding a new node directly below every node of $T$ lying on either the $0$th level or a limit level.
\end{definition}

The following properties are immediate, making use of \cref{res:R-grading iff succ-special}.

\begin{lemma}\label{res:props of shifttree}\
	\begin{enumerate}[label=(\arabic*)]
		\item\label{item:shift succ; res:props of shifttree}
			$\Pad(T)(\suc) = T$.
		\item\label{item:R gradable; res:props of shifttree}
			$\Pad(T)$ is \R-gradable if and only if $T$ is \Q-gradable.
		\item\label{item:uncountable branches; res:props of shifttree}
			$T$ has no uncountable branches if and only if $\Pad(T)$ has no uncountable branches.
		\item\label{item:Hausdorff; res:props of shifttree}
			$T$ is Hausdorff if and only if $\Pad(T)$ is Hausdorff.
	\end{enumerate}
\end{lemma}

The task is thus to construct a Hausdorff, \Q-ungradable well-stratified tree with no uncountable branches. This is achieved by the following definition and result, due to Laver (unpublished; reported in \cite[Ch.~4, Theorem 4(a)]{baumthesis}).

\begin{definition}
	Let $\In_\omega$ be the tree of all injective functions of the form $f \colon \alpha \to \omega$, where $\alpha$ is an ordinal, ordered by $\sse$. That is:
	\begin{equation*}
		\In_\omega = \{f \colon \alpha \to \omega \mid \alpha\text{ is an ordinal and }f\text{ is injective}\}
	\end{equation*}
\end{definition}

\begin{theorem}\label{res:In omega Q-ungradable}
	$\In_\omega$ has no \Q-grading.
\end{theorem}

\begin{proof}
	Following \cref{res:Q-grading iff special}, assume for a contradiction that $\In_\omega = \bigcup_{n \in \omega\setminus \{0\}} A_n$, where each $A_n$ is an antichain. We will construct by induction a sequence of elements $f_0 \subset f_1 \subset \cdots$ of $\In_\omega$, all with coinfinite range, together with a sequence of natural numbers $x_1, x_2, \ldots$ such that $\ran(f_n) \cap \{x_1, \ldots, x_n\} = \es$. Each $x_i$ represents a `promise' that it will never appear in the range of an $f_n$; together they ensure that $\bigcup_{n \in \omega}f_n$ has coinfinite range. 

	Start with $f_0 \defeq \es$. Assume that $f_{n-1}$ is constructed. Choose any $f_{n} \in \In_\omega$ and $x_{n}$ subject to the following conditions.
	\begin{enumerate}[label=(\roman*)]
		\item\label{item:proper; pr:In omega Q-ungradable}
			$f_{n}$ is a proper extension of $f_{n-1}$.
		\item\label{item:coinfinite; pr:In omega Q-ungradable}
			$f_{n}$ has coinfinite range.
		\item\label{item:cumplir; pr:In omega Q-ungradable}
			$\ran(f_{n}) \cap \{x_1, \ldots, x_{n}\} = \es$.
		\item\label{item:in A if poss; pr:In omega Q-ungradable}
			Choose $f_n \in A_n$ if this is possible for some function and $x_n$ satisfying conditions \ref{item:proper; pr:In omega Q-ungradable}, \ref{item:coinfinite; pr:In omega Q-ungradable} and \ref{item:cumplir; pr:In omega Q-ungradable}.
	\end{enumerate}
	Note that \ref{item:proper; pr:In omega Q-ungradable}, \ref{item:coinfinite; pr:In omega Q-ungradable} and \ref{item:cumplir; pr:In omega Q-ungradable} can be satisfied since $f_{n-1}$ has coinfinite range.

	Now let $f \defeq \bigcup_{n \in \omega} f_n$. Then $f \in \In_\omega$, so $f \in A_n$ for some $n > 0$. But note that $\ran(f) \cap \{x_1, x_2, \ldots\} = \es$, so $f$ satisfies \ref{item:proper; pr:In omega Q-ungradable}, \ref{item:coinfinite; pr:In omega Q-ungradable} and \ref{item:cumplir; pr:In omega Q-ungradable} above at stage $n$. Hence by \ref{item:in A if poss; pr:In omega Q-ungradable} we must have that $f_n \in A_n$. Then $f_n \subset f$ contradicts that $A_n$ is an antichain. \contradiction
\end{proof}

Putting it all together, we can prove \cref{res:non cont grad brto}, thus answering the main question in the negative.

\begin{proof}[Proof of \cref{res:non cont grad brto}]
	Take the tree $\Pad(\In_\omega)$. By \cref{res:In omega Q-ungradable,res:props of shifttree}, this is an \R-ungradable well-stratified tree with no uncountable branches. To see that it is Hausdorff, by \cref{res:props of shifttree}\ref{item:Hausdorff; res:props of shifttree} it suffices to show that $\In_\omega$ is Hausdorff. Take $f, g \in \In_\omega$ distinct. Let $\alpha$ be the least ordinal at which $f$ and $g$ disagree. Then $f \restr \alpha$ is the meet of $f$ and $g$.

	Therefore, by \cref{res:road space of cb tree is brot}, the road space, $\Road(\Pad(\In_\omega))$ is a branchwise-real tree order. Furthermore, it has no continuous grading, since any such grading would restrict to an \R-grading on $\Pad(\In_\omega)$, via the canonical embedding.
\end{proof}

\Cref{res:brot param iff every well-strat grad} below shows that the use of an \R-ungradable well-stratified tree with no uncountable branches is in some sense essential: every branchwise-real tree order with no continuous grading contains such a well-stratified tree.

\section{The connection to well-stratified trees}
\label{sec:inside}

In this section, I deepen the connection between branchwise-real tree orders and well-stratified trees with the following theorem.

\begin{theorem}\label{res:brot param iff every well-strat grad}
	A branchwise-real tree order $X$ is continuously gradable if and only if every well-stratified subtree $T \sse X$ is \R-gradable.
\end{theorem}

The proof of the non-trivial direction works by producing a sequence of well-stratified subtrees approximating $X$. The following is the main construction.

\begin{construction}\label{con:sequence of subtrees}
	Let $X$ be a branchwise-real tree order with root $p$. Fix a well-ordering $(B_\alpha \mid \alpha < \kappa)$ of the branches of $X$. For each $\alpha$ let $F_\alpha \defeq B_\alpha \setminus \bigcup_{\beta < \alpha} B_\beta$ be the final segment of $B_\alpha$ disjoint from the previous branches. Since $X$ is a branchwise-real tree order, each $F_\alpha$ is isomorphic to a real interval (which may be empty or a singleton). Fix an isomorphism $r_\alpha \colon F_\alpha \to I_\alpha$, such that $I_\alpha$ is either empty, the singleton $\{1\}$ or a unit interval. Fix an enumeration $(q_n \mid n \in \omega)$ of $\Q \cap [0,1]$. Let $T_n$ be the root $p$ of $T$ together with the union of the preimages according to each $r_\alpha$ of the set $\{q_0, \ldots, q_n\}$:
	\begin{equation*}
		T_n \defeq \{p\} \cup \bigcup_{\alpha < \kappa}r_\alpha^{-1}\left(\{q_0, \ldots, q_n\}\right)
	\end{equation*}
\end{construction}

\begin{lemma}\label{res:T_ns bc and well-strat}
	Each $T_n$ in \cref{con:sequence of subtrees} is a well-stratified tree with no uncountable branches.
\end{lemma}

\begin{proof}
	Note that the segments $F_\alpha$ partition $X$. Take any $x \in T_n$. Let us see directly that $\ds x$ is well-ordered. Take any non-empty $S \sse \ds x$. Note that, since $S$ is linearly ordered, if $y \in F_\alpha \cap S$ and $z \in F_\beta \cap S$ with $\alpha < \beta$ then $y < z$. Let $\alpha < \kappa$ be least such that $F_\alpha \cap S \neq \es$. Since $F_\alpha \cap T_n$ is finite, the set $F_\alpha \cap S$ has a least element, which is then least in $S$.

	To see that $T_n$ has no uncountable branches, note that each branch $B$ of $T_n$ is a subset of a branch in $X$, and thus embeds into \R. Since no uncountable ordinal order-embeds into \R, we must have that $B$ is countable.
\end{proof}

Because we eventually add every element in the preimages of $\Q$ to the $T_n$'s, the union $\bigcup_{n \in \omega}T_n$ will be a dense subtree of $X$ with respect to the interval topology, as follows.

\begin{definition}
	Let $X$ be a tree order with root $p$. The \emph{interval topology} on $X$ is the topology generated by taking the intervals $[p,x)$ for $x \in X$, $(x,y)$ for $x < y$ in $X$ and $(x,y]$ for $x<y$ in $X$ with $y$ a maximum element, as the basic open sets.\footnote{Recall that $(x,y)$ is the interval $\{z \in X \mid x < z < y\}$, etc.}
\end{definition}

\begin{lemma}\label{res:union of T_ns dense}
	Let $X$ be a branchwise-real tree order and let $(T_n)_{n \in \omega}$ be as in \cref{con:sequence of subtrees}. Then $\bigcup_{n \in \omega}T_n \sse X$ is dense in the interval topology on $X$. In other words, for any $x < y$ in $X$ there is $n \in \omega$ and $z \in T_n$ with $x < z < y$.
\end{lemma}

\begin{proof}
	Take $x < y$. Let $\alpha$ be least such that $F_\alpha \cap (x,y) \neq \es$. By the minimality of $\alpha$, the set $F_\alpha \cap (x,y)$ is an initial segment of $(x,y)$, and so is isomorphic to a non-trivial real interval. But then it must contain some $z \in r_\alpha^{-1}(\Q)$. This $z$ eventually appears in some $T_n$, by construction, and satisfies $x < z < y$.
\end{proof}

The last result we need to establish the connection between branchwise-real tree orders and our increasing sequences of well-stratified trees is that any continuous grading of $X$ corresponds to an \R-grading on each of the $T_n$'s. For this, it is first necessary to show that the continuous gradability of $X$ is equivalent to the apparently weaker notion of simple \R-gradability. That is, any \R-grading of a branchwise-real tree order can be transformed into a continuous grading. This is done by `removing all the gaps', as follows.

\begin{theorem}\label{res:para iff R-grading}
	Let $X$ be a branchwise-real tree order. Then $X$ has a continuous grading if and only if it has an \R-grading.
\end{theorem}

\begin{proof}
	The left-to-right is immediate. So assume that $f \colon X \to \R$ is an \R-grading. We will go through the tree eliminating all the discontinuities in $f$, in a Zorn's Lemma style argument. We work with the set of partial monotonic functions $\ell \colon X \partto \R$ such that:
	\begin{enumerate*}[label=(\alph*)]
		\item $\dom(\ell)$ is downwards-closed,
		\item $\ell$ is continuous on its domain, in the sense of \cref{def:continuous grading}, and
		\item $\ell \leq f$ on its domain.
	\end{enumerate*}
	By Zorn's Lemma, there is a maximal such partial function $\ell \colon X \partto \R$. Suppose for a contradiction that the domain of $\ell$ is not $X$. 

	Pick some maximal linearly ordered subset $C \sse X \setminus \dom(\ell)$. If $C$ consists of a single point $x$, we can extend $\ell$ to $\wh\ell \colon \dom(\ell) \cup \{x\} \to \R$ by setting $\wh\ell(x) \defeq \sup_{y < x} \ell(y)$, noting that this is bounded by $f(x)$. Therefore, we may assume that $C$ is not a singleton. Now, since $\dom(\ell)$ is downwards-closed, $C$ is a final segment in some branch $B$ of $X$. Hence there is an isomorphism $r \colon C \to I$ onto a unit interval. The map $f \circ r^{-1} \colon I \to \R$ is then monotonic. By a result coming from real analysis, the only discontinuities on $f \circ r^{-1}$ are jump discontinuities, and there are at most countably many such (see for instance Theorem 4.30 of \cite{baby}). We will find some continuous monotonic function lying below $f \circ r^{-1}$. See \cref{fig:f on C; proof:para iff R-grading} for a picture.

	\begin{figure}[ht]
		\begin{equation*}
			\begin{tikzpicture}
				\begin{scope}
					\filldraw[lightgray] (0,0) node[point] {} -- (-2,4) -- (0,4) --cycle; 
					\draw[dotted] (0,0) -- (-2,4);
					\draw[dotted] (0,0) -- (2,4);
					\draw (0,3) node (Ca) {} -- (1.5,4) node (Cb) {};
					\draw[brace]
						($(Cb) + (0.1,-0.1)$) -- node[below right] {$C$} ($(Ca) + (0.1,-0.1)$);
					\node at (-0.75,3) {$\dom(\ell)$};
				\end{scope}
				\begin{scope}[xshift=-120, yshift=-110, scale=1.5]
					\draw[->] (0,0) -- node[left] {$f$} (0,2);
					\draw (0,0) node (a1a) {} -- (2,0) node (a1b) {};
					\draw[brace]
						($(a1b) + (0,-0.1)$) -- node[yshift=-12] {$C$} ($(a1a) + (0,-0.1)$);
					\begin{scope}
						\draw (0,0) node (l1a) {} .. controls +(0.2,0) and (0.3,0.2) .. (0.5,0.2) node (l1b) {};
					\end{scope}
					\begin{scope}[yshift=10]
						\draw (0.5,0.2) node (l1c) {} .. controls +(0.3,0) .. (1,0.3) node (l1d) {};
					\end{scope}
					\begin{scope}[yshift=15]
						\draw (1,0.3) node (l1e) {} .. controls +(0,0.3) .. (1.5,0.6) node (l1f) {};
					\end{scope}
					\begin{scope}[yshift=30]
						\draw (1.5,0.6) node (l1g) {} -- (2,0.8) node (l1h) {};
					\end{scope}
					\begin{scope}[every path/.style={lightgray, dotted}]
						\draw (l1b) -- (l1c);
						\draw (l1d) -- (l1e);
						\draw (l1f) -- (l1g);
					\end{scope}
				\end{scope}
				\begin{scope}[xshift=40, yshift=-110, scale=1.5]
					\draw[->] (0,0) -- node[left] {$\wh \ell$} (0,2);
					\draw (0,0) node (a2a) {} -- (2,0) node (a2b) {};
					\draw[brace]
						($(a2b) + (0,-0.1)$) -- node[yshift=-12] {$C$} ($(a2a) + (0,-0.1)$);
					\begin{scope}[yslant=0.35]
						\draw (0,0) node (l1a) {} .. controls +(0.2,0) and (0.3,0.2) .. (0.5,0.2) node (l1b) {};
						\draw (0.5,0.2) node (l1c) {} .. controls +(0.3,0) .. (1,0.3) node (l1d) {};
						\draw (1,0.3) node (l1e) {} .. controls +(0,0.3) .. (1.5,0.6) node (l1f) {};
						\draw (1.5,0.6) node (l1g) {} -- (2,0.8) node (l1g) {};
					\end{scope}
				\end{scope}
				\begin{scope}[xshift=0, yshift=-110, scale=1.5]
					\draw[decorate, decoration=snake, ->] (-0.4,0.8) -- (0.4,0.8);
				\end{scope}
			\end{tikzpicture}
		\end{equation*}
		\caption{Extending $\ell$ to $\wh\ell \colon \dom(\ell) \cup C \to \R$ by modifying $f$ on $C$}
		\label{fig:f on C; proof:para iff R-grading}
	\end{figure}
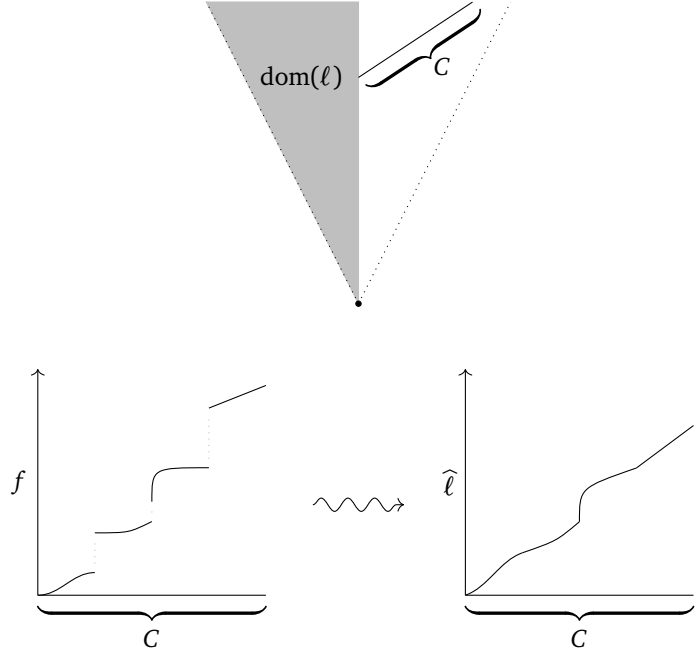

	For $x \in C$, in analogy with real analysis define $f(x-) \defeq \sup_{y < x} f(y)$. Further, define the \emph{jump at $x$ on $C$} as:
	\begin{equation*}
		j_C(x) \defeq \inf_{\substack{y \in C\\ y > x}} f(y) - f(x-)
	\end{equation*}
	By the above, there are at most countably many $x$'s on $C$ such that $j_C(x) > 0$. To extend $\ell$ to $C$, we first remove each such discontinuity, defining $g_C \colon C \to \R$ by:
	\begin{equation*}
		g_C(x) \defeq f(x-) - \sum_{\substack{y \in C\\ y < x}} j_C(y) 
	\end{equation*}
	Then $g_C$ is continuous with respect to $r$ (i.e.\@ $g_C \circ r^{-1} \colon I \to \R$ is continuous). It is also weakly monotonic, but could fail to be strictly monotonic, e.g.\@ in the case where the jump discontinuities are dense. To obtain a strictly monotonic function, we add a contribution from each jump in a continuous way. Define $h_C \colon C \to \R$ by:
	\begin{equation*}
		h_C(x) \defeq g_C(x) + \sum_{\substack{y \in C\\ y < x}}\frac{r(x)-r(y)}{1-r(y)}j_C(y)
	\end{equation*}
	Then function $h_C \circ r^{-1}$ is then the uniform limit of continuous functions, and therefore continuous. Note also that $h_C \leq f$ on $C$.

	Finally, it remains to extend $\ell$ to $C$ by attaching $h_C$. For this we may need to shift $h_C$  a little so that it fits in continuously with $\ell$. Extend $\ell$ to $\wh \ell \colon \dom(\ell) \cup C \to \R$ by letting, for $x \in C$:
	\begin{equation*}
		\wh \ell(x) \defeq h_C(x) + \left(\sup_{z \in B \setminus C} \ell(z) - \inf_{y \in C} h_C(y)\right)
	\end{equation*}
	The resulting function is then satisfies the continuity condition, and is moreover such that $\wh \ell \leq f$ on its domain. This contradicts the maximality of $\ell$. \contradiction
\end{proof}

Finally, we can establish the connection between the continuous gradability of $X$ and the \R-gradability of each $T_n$ in \cref{con:sequence of subtrees}.

\begin{theorem}\label{res:brot para iff every T_n gradable}
	Let $X$ be a branchwise-real tree order and let $(T_n)_{n \in \omega}$ be as in \cref{con:sequence of subtrees}. Then $X$ is continuously gradable if and only if every $T_n$ is \R-gradable.
\end{theorem}

\begin{proof}
	Firstly, any continuous grading of $X$ restricts to an \R-grading on each $T_n$. Conversely, assume that each $T_n$ has an \R-grading $f_n \colon T_n \to [0,1)$. Now, each function $f_n$ can be extended to a weakly monotonic function $\wh f_n \colon X \to [0,1]$ by:
	\begin{equation*}
		\wh f_n(x) \defeq \sup \{f_n(y) \mid y \in T_n \text{ and }y \leq x\}
	\end{equation*}
	Then define $f \colon X \to [0,1]$ by:
	\begin{equation*}
		f(x) \defeq \sum_{n \in \omega}^\infty \frac{\wh f_n(x)}{2^n}
	\end{equation*}
	Let us see that $f$ it is strictly monotonic. Take $x < y$ in $X$. By \cref{res:union of T_ns dense} there is $n \in \omega$ and $x_0, y_0 \in T_n$ such that $x < x_0 < y_0 < y$. Since $f_n$ is monotonic on $T_n$, we have $\wh f_n(x_0) < \wh f_n(y_0)$. Since also $\wh f_m(x_0) \leq \wh f_m(y_0)$ for all $m \in \omega$, we get that $f(x_0) < f(y_0)$. Hence:
	\begin{equation*}
		f(x) \leq f(x_0) < f(y_0) \leq f(y)
	\end{equation*}
	Therefore $f \colon X \to \R$ is an \R-grading, and thus, by \cref{res:para iff R-grading}, $X$ is continuously gradable.
\end{proof}

This last piece allows us to finish the proof of this section's main result.

\begin{proof}[Proof of \cref{res:brot param iff every well-strat grad}]
	If $\ell \colon X \to \R$ is a continuous grading, then it restricts to an \R-grading of any well-stratified subtree. Conversely, assume that $X$ has no continuous grading. Let $(T_n)_{n \in \omega}$ be as in \cref{con:sequence of subtrees}. Then by \cref{res:brot para iff every T_n gradable}, at least one $T_n$ must be \R-ungradable (in fact, infinitely many are).
\end{proof}

\section{Separability and wispiness}
\label{sec:branching nodes}

With the result of \cref{sec:inside} established, we can now move on to consider refinements of the main question. Our original question was motivated by considering the underlying orders of \R-trees. What happens if we look instead at \emph{separable} \R-trees? \Cref{res:seperability for brots} below shows that an \R-tree is separable if and only if any cut-point order is what I call `countably wispy': that it contains fewer than countably many branching nodes, and fewer than countably many terminal segments isomorphic to an non-trivial interval. The corresponding refinement of the main question then becomes the following. In contrast to the main question, the answer turns out to be `Yes'.

\begin{question}\label{q:countably wispy}
	Is every countably wispy branchwise-real tree order continuously gradable?
\end{question}

The first item of business is to make the term `countably wispy' precise, and establish its correspondence with separability.

\begin{definition}
	Let $X$ be a tree order and $x \in X$. An \emph{$x$-connected component} is an equivalence class of $\us x \setminus \{x\}$ under the relation:
	\begin{equation*}
		y \sim_x z \quad\Lra\quad \text{there is } w > x \text{ such that } w \leq y,z
	\end{equation*}
	The \emph{degree}, $\deg(x)$, of $x$ is the number of $x$-connected components. Say that $x$ is \emph{terminal} if $\deg(x) = 0$. Say that $x$ is \emph{branching} if $\deg(x)>1$. Define $\Branch(X)$ to be the set of all branching nodes in $X$.
\end{definition}

Note that we can endow any partial order with the so-called `Alexandrov topology', in which the open sets are exactly those which are upwards-closed. Under this topology, the $x$-connected components are precisely the maximal connected subsets of $\us x \setminus \{x\}$.

\begin{definition}
	A \emph{twig} in a branchwise-real tree order is a maximal upwards-closed subset isomorphic to a non-trivial real interval.
\end{definition}

\begin{figure}[t]
	\begin{equation*}
		\begin{tikzpicture}[scale=0.8]
			\filldraw[lightgray] (0,0) -- (-2,4) -- (2,4) -- cycle;
			\node[point] at (0,0) [label=right:{root}] {};
			\node (rt) at (3,2) {rest of the tree};
			\draw[->] (rt) -- (0.5,3);
			\draw (-1,2) -- (-3,4);
			\draw[brace] 
				($(-1,2)+(-0.2,-0.2)$) -- node[below left] {twig} ($(-3,4)+(-0.2,-0.2)$);
		\end{tikzpicture}
	\end{equation*}
	\caption{What a twig looks like}
	\label{fig:twig}
\end{figure}
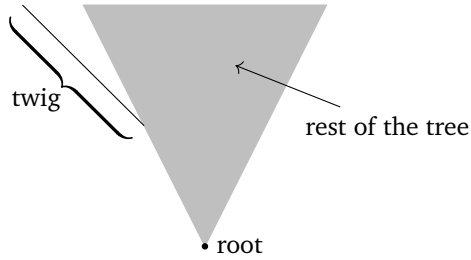

See \cref{fig:twig} for a picture of a twig.

\begin{definition}
	A branchwise-real tree order is a \emph{countably wispy} if it contains at most countably many branching nodes and at most countably many twigs.
\end{definition}

Note that a countably wispy tree may still contain uncountably many terminal nodes. For example, the road space of the well-stratified tree $2^{\leq\omega}$ has continuum-many terminal nodes, but no twigs and only countably many branching nodes.

\begin{lemma}\label{res:seperability for brots}
	An \R-tree $\ab{X,d}$ is separable if and only if each cut-point order is countably wispy.
\end{lemma}

\begin{proof}
	Fix a cut-point order on $X$ with root $p$. Assume $D \sse X$ is countable and dense in the metric topology. First note that for any $x \in \Branch(X)$, any $x$-connected component $C$ and any $y \in C$ the open ball about $y$ of radius $d(x,y)$ is open and contained in $C$; hence $C \cap D \neq \es$. In particular, every branching node sits below an element of $D$. Now, take any $y \in D$. For every branching node $x < y$, by the above method we can pick $x^* \in D$ which is not in the same $x$-connected component as $y$. All the $x^*$'s must be distinct, and hence the number of branching nodes below $y$ must be countable. Thus, as the countable union of countable sets is countable, there can only be countably many branching nodes in total. Finally, note that any two distinct twigs are disjoint, and that any twig must intersect with $D$; hence there are only countably many twigs.

	Conversely, assume that the cut-point order with root $p$ is countably wispy. Let $X_0$ be the set of terminal nodes in $X$, and let $X_{\mr t}$ be the union of the twigs in $X$. For any $y \in X \setminus (X_0 \cup X_{\mr t})$, there is $x \in \Branch(X)$ such that $x \geq y$. Pick for each $x \in \Branch(X)$ a branch $B$ such that $x \in B$. Furthermore, for each twig there is a unique branch which contains it. Thus we can find a countable set $\mc B$ of branches through $X$ such that $X = \bigcup \mc B \cup X_0$. In particular, $\bigcup \mc B$ is dense in $X$. Therefore the set:
	\begin{equation*}
		D \defeq \bigcup_{B \in \mc B} \{y \in B \mid d(p,y) \in \Q\}
	\end{equation*}
	is countable and dense in $X$.
\end{proof}

With the connection thus established, we can now answer \cref{q:countably wispy}. In fact, we only need half of the definition of `countably wispy'.

\begin{theorem}\label{res:BROT BrNo countable param}
	Any branchwise-real tree order containing at most countably many branching nodes is continuously gradable.
\end{theorem}

The proof of this result goes via \cref{res:brot param iff every well-strat grad}: to show that a branchwise-real tree order is continuously gradable, it suffices to show than any well-stratified subtree is \R-gradable. For this, we will make use of the following three general lemmas.

\begin{lemma}\label{res:BrNo well-strat Q-grad}
	If $T$ is a well-stratified tree with no uncountable branches and $\Branch(T)$ is \Q-gradable, then $T$ is \R-gradable.
\end{lemma}

\begin{proof}
	By \cref{res:Q-grading iff special} we have $\Branch(T) = \bigcup_{n \in \omega} A_n$, where each $A_n$ is an antichain. Define $f_0 \colon \Branch(T) \to \R$ by:
	\begin{equation*}
		f_0(x) \defeq \sum \left\{\frac 1 {2^n} \;\middle|\; \exists y \leq x \colon y \in A_n\right\}
	\end{equation*}
	Note that $f_0$ is a bounded \R-grading on $\Branch(T)$ with the property that for any $x \in \Branch(T)$:
	\begin{equation*}
		f_0(x) > \sup_{y<x} f_0(y)
	\end{equation*}
	In other words, below each $x$ there is a `jump' of size $f_0(x) - \sup_{y<x} f_0(y)$. We now extend $f_0$ to $f \colon T \to \R$. Define an equivalence relation on $T \setminus \Branch(T)$, setting $u \sim v$ if and only if either $u \leq v$ or $v \leq u$, and:
	\begin{equation*}
		\ds u \cap \Branch(T) = \ds v \cap \Branch(T)
	\end{equation*}
	Since $T$ has no uncountable branches, each equivalence class is countable; it is moreover well-ordered, and so order-isomorphic to a countable ordinal. Let $C$ be any such class. If there is $x \in \Branch(T)$ such that $x > u$ for all $u \in C$, then there is a least such $x$. By \cref{res:countable ordinal in Q} we can define $f$ monotonic on $C$ such that for all $u \in C$:
	\begin{equation*}
		\sup_{y<x} f_0(y) < f(u) < f_0(x)
	\end{equation*}
	If there is no such $x$, then as $f_0$ is bounded, we can again define $f$ monotonic on $C$ such that for all $u \in C$:
	\begin{equation*}
		\sup_{y<x} f_0(y) < f(u)
	\end{equation*}
	Putting it all together, we arrive are our desired \R-grading $f \colon T \to \R$.
\end{proof}

\begin{lemma}
	If $X$ is a branchwise-real tree order, then $X$ is a complete meet-semilattice: any non-empty $S \sse X$ has a greatest lower bound $\bigwedge S$.
\end{lemma}

\begin{proof}
	Take any $x \in S$ and consider the set $R \defeq \{x \wedge y \mid y \in S\}$. Since $\ds x$ is linearly ordered, so is $R$. Hence $R$ lies in some branch $B$ of $X$. Since $B$ is order-isomorphic to a real interval, $R$ has an infimum in $B$, which infimum is then also the element $\bigwedge S$.
\end{proof}

\begin{lemma}\label{res:BrNo aprox leq BrNo brot}
	Let $X$ be a branchwise-real tree order and let $T \sse X$ be a well-stratified subtree. Then $\abs{\Branch(T)} \leq \abs{\Branch(X)}$.
\end{lemma}

\begin{proof}
	Define a function $s \colon \Branch(T) \to \Branch(X)$ as follows. For any $x \in \Branch(T)$, let $S_x$ be its set of immediate successors in $T$ (of which there are at least $2$). Then let $s(x) \defeq \bigwedge S_x$ in $X$. As $\abs{S_x} \geq 2$, we have $s(x) \in \Branch(X)$. Let us see that $s$ is injective. Assume that $s(x)=s(y)$ but $x \neq y$. Since $x,y \leq s(x)$ in $X$, which is a tree order, we must have that $x$ and $y$ are comparable, say $x < y$. But then $x$ has an immediate successor $z$ in $T$ such that $z \leq y$, and $s(x) < z$ (since $x$ has more than one immediate successor). This means that:
	\begin{equation*}
		s(x) < z \leq y \leq s(y)
	\end{equation*}
	which is a contradiction. \contradiction
\end{proof}

We can now prove the theorem, answering \cref{q:countably wispy}.

\begin{proof}[Proof of \cref{res:BROT BrNo countable param}]
	Let $X$ be branchwise-real tree order containing at most countably many branching nodes. Take any $T \sse X$ a well-stratified tree. By \cref{res:BrNo aprox leq BrNo brot} we have that $\abs{\Branch(T)} \leq \abs{\Branch(X)} \leq \aleph_0$. So $\Branch(T)$ is a countable well-stratified tree. This means that it has countable height $\alpha$. By \cref{res:countable ordinal in Q}, there is a monotonic map $\alpha \to \Q$, which then pulls back through the rank function to a \Q-grading on $\Branch(T)$. Then by \cref{res:BrNo well-strat Q-grad}, the tree $T$ is \R-gradable. Thus, any well-stratified subtree of $X$ is \R-gradable, and so by \cref{res:brot param iff every well-strat grad} we get that $X$ is continuously gradable.
\end{proof}

The notion of countable-wispiness generalises naturally to that of $\kappa$-wispiness, about which we can ask analogous questions to the above.

\begin{definition}
	A branchwise-real tree order is a \emph{$\kappa$-wispy} if it contains fewer than $\kappa$-many branching nodes and fewer than $\kappa$-many twigs.
\end{definition}

On the metric space side, for uncountable $\kappa$ this property corresponds to the following generalisation of separability.

\begin{definition}
	Let $\kappa$ be an infinite cardinal. A topological space is \emph{$\kappa$-separable} if it admits a dense subset of size less than $\kappa$.
\end{definition}

\begin{lemma}\label{res:kappa sep iff kappa wispy}
	Let $\kappa$ be an infinite cardinal. An \R-tree $\ab{X,d}$ is $\kappa$-separable if and only if each cut-point order is $\kappa$-wispy.
\end{lemma}

\begin{proof}
	The only trouble encountered when generalising the proof of \cref{res:seperability for brots} is that, in showing that separability implies that there are only countably many branching nodes, we used that the countable union of countable sets is countable. So we need to be a little careful when $\kappa$ is not regular. But note that if $X$ is $\kappa$-separable for $\kappa$ a limit cardinal, then it is $\lam$-separable for some regular $\lam < \kappa$.
\end{proof}

Now, in contrast with the countable case, the existence of an $\aleph_\alpha$-wispy branchwise-real tree order with no continuous grading, for $\alpha \geq 2$, is independent of \ZFC.

\begin{theorem}\label{res:CH non-para BrNo aleph_1}
	Let $X = \Road(\Pad(\In_\omega))$ be the tree constructed in \cref{sec:non-para}. Then $X$ is $\cont^+$-wispy. Hence under the Continuum Hypothesis there is a $\aleph_2$-wispy branchwise-real tree order with no continuous grading.
\end{theorem}

\begin{proof}
	The structure $\Branch(X)$ is isomorphic to $\Branch(\In_\omega)$. An element of $\In_\omega$ is an injective function from a countable ordinal into $\omega$. There are $\aleph_1$-many countable ordinals, and each has at most $\cont$-many injective functions into $\omega$. Hence:
	\begin{equation*}
		\abs{\Branch(X)} = \abs{\Branch(\In_\omega)} \leq \abs{\In_\omega} = \aleph_1 \cdot \cont = \cont
	\end{equation*}
	Furthermore, a twig in $X$ corresponds to a terminal node in $\In_\omega$, of which there are $\cont$-many. Thus $X$ is $\cont^+$-wispy.
\end{proof}

\begin{theorem}\label{res:MA kappa BrNo para}
	Let $\kappa$ be a cardinal. Under $\MA_\kappa$ every branchwise-real tree order $X$ with at most $\kappa$-many branching nodes is continuously gradable.
\end{theorem}

\begin{proof}
	Note that when $\kappa \leq \aleph_0$ this already follows from \cref{res:BROT BrNo countable param}. Take any $T \sse X$ a well-stratified subtree. By \cref{res:brot param iff every well-strat grad}, it suffices to show that $T$ is \R-gradable. By \cref{res:BrNo well-strat Q-grad}, it suffices to find a \Q-grading of $\Branch(T)$. By \cref{res:BrNo aprox leq BrNo brot}, we have that $\abs{\Branch(T)} \leq \abs{\Branch(X)} \leq \kappa$. We find a \Q-grading for $\Branch(T)$ using the usual `specialising forcing'.

	Let $\bb P$ be the partial order consisting of all finite, monotonic partial functions $p \colon \Branch(T) \partto \Q$ under:
	\begin{equation*}
		p \leq q \quad\Lra\quad p\text{ extends }q
	\end{equation*}
	It is well-known that, as $\Branch(T)$ has no uncountable branches, this poset is ccc (see for example Theorem~III.5.19 in \cite{kunen}). Furthermore, for any $x \in \Branch(T)$ the following set is dense in $\bb P$.
	\begin{equation*}
		D_x \defeq \{p \in \bb P \mid x \in \dom(p)\}
	\end{equation*}
	By $\MA_\kappa$ then, using that $\abs{\Branch(T)} \leq \kappa$, there is a filter $G$ on $\bb P$ which intersects with every $D_x$. Letting $f \defeq \bigcup G$, we see that it is a monotonic function $\Branch(T) \to \Q$; in other words, a \Q-grading.
\end{proof}

\begin{remark}
	An alternative way of showing that $X$ is continuously gradable is as follows. Enumerate the branches of $X$ as $B_\alpha \mid \alpha < \lam$, and let $F_\alpha \defeq B_\alpha \setminus \bigcup_{\beta < \alpha} B_\beta$. For each $\alpha$ fix an order isomorphism $r_\alpha \colon F_\alpha \to I_\alpha$ onto a real interval. Then let:
	\begin{equation*}
		X_\Q \defeq \bigcup_{\alpha < \lam} r_\alpha^{-1}(\Q) 
	\end{equation*}
	This is a tree order in which every branch is countable. We can then take the specialising forcing of finite, monotonic partial functions $X_\Q \partto \Q$. It is possible to show that this forcing poset is ccc (e.g.\@ by adapting the proof of Theorem~III.5.19 in \cite{kunen}), and hence by a density argument $\MA_\kappa$ entails that there is a \Q-grading $X_\Q \to \Q$. This can then be completed to an \R-grading of $X$, which by \cref{res:para iff R-grading} means that $X$ must be continuously gradable.
\end{remark}

\begin{corollary}\label{res:kappa-wispy indep}
	Assume that \ZFC\ is consistent. For $\alpha \geq 2$ it is consistent that there exists an $\aleph_\alpha$-wispy branchwise-real tree order with no continuous grading, and it is also consistent that such a tree does not exist.
\end{corollary}

\begin{proof}
	It is a classical result that the Continuum Hypothesis is consistent with $\ZFC$ (see for instance Theorem~13.20 in \cite{jech}). Hence by \cref{res:CH non-para BrNo aleph_1} it is consistent that there exists a $\aleph_2$-wispy branchwise-real tree order with no continuous grading. Such a tree is also $\kappa$-wispy for any $\kappa \geq \aleph_2$. Conversely, $\MA_{\aleph_{\alpha}}$ is consistent with $\ZFC$, and hence by \cref{res:MA kappa BrNo para} so is the statement that all $\aleph_{\alpha}$-wispy branchwise-real tree orders are continuously gradable.
\end{proof}

\section{The countable chain condition}
\label{sec:suslin}

When building our non-continuously-gradable branchwise-real tree order in \cref{sec:non-para}, we looked for an \R-ungradable well-stratified tree with no uncountable branches. There we constructed one explicitly in \ZFC, but an observant reader may have noticed that a Suslin tree also satisfies these requirements. Motivated by this observation, in this section we consider branchwise-real tree orders satisfying the countable chain condition and its generalisations. The principle question is the following refinement of the main question.

\begin{question}\label{q:ccc}
	Is every ccc branchwise-real tree order continuously gradable?
\end{question}

\begin{remark}\label{rem:well-strat ccc R-grad}
	Note that in the well-stratified case, ccc-ness plus \R-gradability implies that the tree has countable height. Indeed, any ccc, \R-gradable well-stratified tree of height $\omega_1$ must be a Suslin tree, but no Suslin tree admits an \R-grading (see \cite[Lemma~IV.6.5]{kunen}).
\end{remark}

Before delving into the details, let us first see what property of \R-trees corresponds to ccc-ness.

\begin{definition}
	Let $X$ be a connected topological space. A \emph{noncut-point} in $X$ is an element $x \in X$ such that $X \setminus \{x\}$ is connected.
\end{definition}

\begin{theorem}\label{res:cut-point ccc char}
	Let $X$ be an \R-tree. The following are equivalent.
	\begin{enumerate}[label=(\arabic*)]
		\item\label{item:sep noncut; res:cut-point ccc char}
			$X$ is separable and has only countably many noncut-points.
		\item\label{item:count wisp term; res:cut-point ccc char}
			Every cut-point order on $X$ is countably wispy and has only countably many terminal nodes.
		\item\label{item:ccc; res:cut-point ccc char}
			Every cut-point order on $X$ is ccc.
	\end{enumerate}
\end{theorem}

For this we make use of the following lemma due to Bowditch, which shows that on \R-trees connectedness coincides with the natural notion of convexity.

\begin{lemma}\label{res:sub R-tree con iff geodesics}
	A subset $Y \sse X$ of an \R-tree is connected if and only if $[x,y] \sse Y$ for every $x,y \in Y$.\footnote{Recall that $[x,y]$ denotes the unique geodesic segment between $x$ and $y$.}
\end{lemma}

\begin{proof}
	This is Lemma~1.4 of \cite{bowditch99}.
\end{proof}

\begin{proof}[Proof of \cref{res:cut-point ccc char}]
	Throughout, we fix $p \in X$ and consider the cut-point order on $X$ with root $p$.

	$\text{\ref{item:sep noncut; res:cut-point ccc char}} \Lra \text{\ref{item:count wisp term; res:cut-point ccc char}}$.
	By \cref{res:seperability for brots}, $X$ is separable if and only if every cut-point order is countably wispy. Let $X_0$ be the set of terminal nodes with respect to the cut-point order. It suffices then to show that the set of noncut-points is $X_0$ plus possibly $p$. By \cref{res:sub R-tree con iff geodesics}, it is clear that every element of $X_0$ is a noncut-point. Conversely, take any noncut-point $x$ in $X$ other than $p$. If $x$ is not a terminal node, then there is $y \in X$ with $p < x < y$. But then, in $X \setminus \{x\}$, we have $p$ lying in a different connected component to $y$, since the geodesic $[p,y]$ is not a subset of $X \setminus \{x\}$. \contradiction

	$\text{\ref{item:ccc; res:cut-point ccc char}} \Ra \text{\ref{item:count wisp term; res:cut-point ccc char}}$.
	That there are only countably many terminal nodes and twigs is immediate. Let $\ell \colon X \to \R$ be the continuous grading on $X$ found using \cref{res:cut-point para brots}. For each $q \in \Q$, the set $\ell^{-1}\{q\}$ is an antichain, which is therefore countable. Below each such antichain there can only countably many branching nodes, and every branching node must lie below such an antichain. Thus in total there are only countably many branching nodes.

	$\text{\ref{item:count wisp term; res:cut-point ccc char}} \Ra \text{\ref{item:ccc; res:cut-point ccc char}}$.
	Take $A \sse X$ an antichain. Let $A'$ be the result of removing every element of $A$ which is terminal or lies on a twig. Since each twig can contain at most one element of $A$, it suffices to show that $A'$ is countable. For each element $x \in A'$, we can pick $x^* \geq x$ which is a branching node. For $x \neq y$ we must have $x^* \neq y^*$, since otherwise $x$ and $y$ would be comparable. Since there are only countably many branching nodes, we have that $A'$ is countable.
\end{proof}

\begin{remark}
	The equivalence $\text{\ref{item:count wisp term; res:cut-point ccc char}} \Lra \text{\ref{item:ccc; res:cut-point ccc char}}$ shows that on continuously gradable branchwise-real tree orders, ccc-ness is only just stronger than countable wispiness. Note however that the proof of the direction $\text{\ref{item:ccc; res:cut-point ccc char}} \Ra \text{\ref{item:count wisp term; res:cut-point ccc char}}$ makes key use of the continuous grading. We cannot therefore use \cref{res:BROT BrNo countable param} to show that every ccc branchwise-real tree order is continuously gradable.
\end{remark}

We turn to \cref{q:ccc}. The following two theorems show that the answer is independent of \ZFC. Specifically, there exists a ccc branchwise-real tree order with no continuous grading if and only if there exists a Suslin tree.

\begin{theorem}\label{res:road Suslin is ccc non-para brot}
	The road space of a Suslin tree is a ccc branchwise-real tree order with no continuous grading.
\end{theorem}

\begin{proof}
	Let $T$ be a Suslin tree. As noted in \cref{rem:well-strat ccc R-grad}, $T$ admits no \R-grading. Therefore $\Road(T)$ has no continuous grading. Furthermore, if $A \sse \Road(T)$ is any antichain, then by shifting the elements of $A$ up a little, we may assume that $A$ is contained in the canonical embedded copy of $T$ in $\Road(T)$ (i.e.\@ $T \times \{1\}$). Therefore, $A$ must be countable, and $\Road(T)$ must be ccc.
\end{proof}

\begin{theorem}\label{res:every ccc non-para brot contains Suslin}
	Every ccc branchwise-real tree order with no continuous grading contains a Suslin tree.
\end{theorem}

\begin{proof}
	Let $X$ be a ccc branchwise-real tree order with no continuous grading. By \cref{res:brot param iff every well-strat grad}, there is $T \sse X$ an \R-ungradable well-stratified subtree with no uncountable branches. Since $T$ is a suborder of $X$, it is also ccc. Furthermore, it must have height $\omega_1$. If not, it would have height $\alpha$ a countable ordinal. Since, by \cref{res:countable ordinal in Q}, $\alpha$ embeds into $\Q$, this would yield a \Q-grading on $T$. \contradiction Thus $T$ is a Suslin tree.
\end{proof}

The countable chain condition generalises to the $\kappa$ chain condition, which leads to a natural generalisation of \cref{q:ccc}.

\begin{definition}
	A tree order $X$ is \emph{$\kappa$-cc}, where $\kappa$ is a cardinal, if $X$ has no antichains of size $\kappa$.
\end{definition}

With a little extra work, we obtain a generalisation of \cref{res:cut-point ccc char} as follows.

\begin{theorem}\label{res:cut-point kappa-cc char}
	Let $X$ be an \R-tree, and let $\kappa$ be an uncountable cardinal. The following are equivalent.
	\begin{enumerate}[label=(\arabic*)]
		\item\label{item:sep noncut; res:cut-point kappa-cc char}
			$X$ is $\kappa$-separable and has fewer than $\kappa$-many noncut-points.
		\item\label{item:count wisp term; res:cut-point kappa-cc char}
			Every cut-point order on $X$ is $\kappa$-wispy and has fewer than $\kappa$-many terminal nodes.
		\item\label{item:ccc; res:cut-point kappa-cc char}
			Every cut-point order on $X$ is $\kappa$-cc.
	\end{enumerate}
\end{theorem}

\begin{proof}
	All directions in the proof of \cref{res:cut-point ccc char} readily generalise, except for $\text{\ref{item:ccc; res:cut-point ccc char}} \Ra \text{\ref{item:count wisp term; res:cut-point ccc char}}$, in which we use that the union of the countably many branching nodes below each $\ell^{-1}\{q\}$ is countable. We therefore need a new argument in the case where $\kappa$ is of countable cofinality. I will show that for such a $\kappa$, any $\kappa$-cc, continuously gradable branchwise-real tree order has fewer than $\kappa$-many branching nodes.

	So, let $\kappa$ be uncountable with $\cf(\kappa) = \aleph_0$. Then, we can choose a sequence $\kappa_1 < \kappa_2 < \cdots$ of uncountable regular cardinals with limit $\kappa$. Let $X$ be a continuously gradable branchwise-real tree order which contains at least $\kappa$-many branching nodes. We need to find an antichain in $X$ of size $\kappa$. This will be constructed in stages $A_0, A_1, \ldots$, all of which are maximal antichains, such that $\abs{A_n} \geq \kappa_n$ for $n \geq 1$.

	Start with $A_0$ any maximal antichain. Now assume that we have constructed $A_n$. If $\abs{A_n} \geq \kappa$ we can stop the construction here. Otherwise, consider the following partition of $X$:
	\begin{equation*}
		X = ((\dS A_n) \setminus A_n) \cup \bigcup_{x \in A_n} \us x
	\end{equation*}
	The number of branching nodes in $(\dS A_n) \setminus A_n$ is at most $\abs {A_n}$. Since the number of branching nodes in $X$ is at least $\kappa$, and $\abs{A_n} < \kappa$, there must be $x_n \in A_n$ such that the number of branching nodes in $\us {x_n}$ is at least $\kappa_{n+1}$. Now, $\us{x_n}$ is a continuously gradable branchwise-real tree order, and so by the generalisation of the proof of $\text{\ref{item:ccc; res:cut-point ccc char}} \Ra \text{\ref{item:count wisp term; res:cut-point ccc char}}$ in \cref{res:cut-point ccc char} to the regular cardinal $\kappa_{n+1}$, there must be a maximal antichain $A_{n+1}' \sse \us{x_n}$ of size at least $\kappa_{n+1}$. Let:
	\begin{equation*}
		A_{n+1} \defeq A_n \cup A_{n+1}' \setminus \{x_n\}
	\end{equation*}
	Finally, let:
	\begin{equation*}
		A \defeq \bigcup_{n \in \omega}A_n \setminus \{x_0, x_1, \ldots\}
	\end{equation*}
	This is then an antichain in $X$ of size $\kappa$, as required.
\end{proof}

\begin{remark}
	As a corollary of this result, we get that if $X$ is $\kappa$-cc for $\kappa$ some singular cardinal, then it is $\lam$-cc for some $\lam < \kappa$. This is a specific case of a more general result due to Paul Erdös and Alfred Tarski \cite{10.2307/1968767}.
\end{remark}

To show that the existence of a $\kappa$-cc branchwise-real tree order with no continuous grading is independent of \ZFC, we can consider the following notion, which coincides with the notion of Suslin tree for $\kappa=\aleph_1$.

\begin{definition}
	Let $\kappa$ be an uncountable cardinal. A \emph{${<}\kappa$-wide \R-ungradable tree} is an \R-ungradable, $\kappa$-cc well-stratified tree with no uncountable branches.
\end{definition}

\begin{proposition}\label{res:ZFC facts about kappa wide Suslin}\
	\begin{enumerate}[label=(\arabic*)]
		\item\label{item:Suslin; res:ZFC facts about kappa wide Suslin}
			A Suslin tree is exactly a ${<}\aleph_1$-wide $\R$-ungradable tree.
		\item\label{item:MA_kappa; res:ZFC facts about kappa wide Suslin}
			$\MA_\kappa$ implies there are no ${<}\kappa^+$-wide $\R$-ungradable trees, for $\kappa \geq \aleph_1$.
		\item\label{item:MA; res:ZFC facts about kappa wide Suslin}
			If $\MA + \neg\CH$ holds, then there are no ${<}\cont$-wide $\R$-ungradable trees.
		\item\label{item:cont plus; res:ZFC facts about kappa wide Suslin}
			There is a ${<}\cont^+$-wide $\R$-ungradable tree.
	\end{enumerate}
\end{proposition}

\begin{proof}
	\begin{enumerate}[label=(\arabic*)]
		\item This follows since, as noted in of \cref{rem:well-strat ccc R-grad}, all Suslin trees are \R-ungradable.

		\item Let $T$ be any $\kappa^+$-cc well-stratified tree of height $\omega_1$ with no uncountable branches. We use the specialising forcing from \cref{res:MA kappa BrNo para} to show that $T$ has a \Q-grading (and hence an \R-grading). Let $\bb P$ be the partial order consisting of all finite, monotonic partial functions $p \colon T \partto \Q$ under:
		\begin{equation*}
			p \leq q \quad\Lra\quad p\text{ extends }q
		\end{equation*}
		As in the proof of \cref{res:MA kappa BrNo para}, $\bb P$ is ccc, and for each $x \in T$ the following set is dense in $\bb P$.
		\begin{equation*}
			D_x \defeq \{p \in \bb P \mid x \in \dom(p)\}
		\end{equation*}
		To count the size of $T$, note that each level $T(\alpha)$ has size at most $\kappa$, and that there are $\aleph_1$-many levels; hence $\abs T \leq \kappa \cdot \aleph_1 = \kappa$. Therefore, by $\MA_\kappa$, there is a filter $G$ which intersects $D_x$ for every $x \in T$. Then $\bigcup G$ is a \Q-grading $T \to \Q$.

		\item Let $T$ be any $\cont$-cc well-stratified tree of height $\omega_1$ with no uncountable branches. Then $T$ is the $\aleph_1$-union of its levels, each of which have size less than $\cont$. Now, \MA\ implies that $\cont$ is regular (see e.g.\@ \cite[p.~176]{kunen}). In particular it has cofinality strictly greater than $\aleph_1$, and so the size of $T$ must be less than $\cont$. Since \MA\ holds, by using the specialising forcing on $T$ we obtain a \Q-grading $T \to \Q$.

		\item Let $T = \Pad(\In_\omega)$ be the well-stratified tree constructed in \cref{sec:non-para}. By the proof of \cref{res:non cont grad brto}, $T$ is \R-ungradable and has no uncountable branches. As in \cref{res:CH non-para BrNo aleph_1}, $T$ has size $\cont$, and hence it is $\cont^+$-cc.\qedhere
	\end{enumerate}
\end{proof}

\begin{corollary}\label{res:kappa wide Suslin independent}
	Assume that \ZFC\ is consistent. For any $\alpha\geq 1$, it is consistent that there exists a ${<}\aleph_\alpha$-wide $\R$-ungradable tree, and it is also consistent that such a tree does not exist.
\end{corollary}

\begin{proof}
	It is consistent that a Suslin tree exists. If this is the case, then by \cref{res:ZFC facts about kappa wide Suslin} \ref{item:Suslin; res:ZFC facts about kappa wide Suslin} there is a ${<}\aleph_1$-wide $\R$-ungradable tree, which is also ${<}\kappa$-wide for all $\kappa \geq \aleph_1$. Conversely, it is consistent that $\MA_{\aleph_{\alpha}}$ holds, and in this case by \cref{res:ZFC facts about kappa wide Suslin} \ref{item:MA_kappa; res:ZFC facts about kappa wide Suslin} there are no ${<}\aleph_\alpha^+$-wide $\R$-ungradable trees. Since any ${<}\aleph_\alpha$-wide $\R$-ungradable tree is a ${<}\aleph_\alpha^+$-wide $\R$-ungradable tree, this implies that the former does not exist either.
\end{proof}

Finally, the proofs of \cref{res:road Suslin is ccc non-para brot,res:every ccc non-para brot contains Suslin} readily generalise, yielding the following independence result.

\begin{theorem}\label{res:bound-kappa to well-strat}
	Let $\kappa$ be an uncountable cardinal. There exists a $\kappa$-cc branchwise-real tree order with no continuous grading if and only if there exists a ${<}\kappa$-wide $\R$-ungradable tree.
\end{theorem}

\section{Open questions}
\label{sec:conclusion}

Another class of trees finding many applications in geometric group theory is that of \emph{\Lam-trees}, a generalisation of \R-trees (see \cite{Shalen1987,Shalen91,Morgan1992LambdaA,ChiswellIan2001Itl}). Briefly, an \emph{ordered abelian group} is an abelian group $\ab{\Lam,+}$ equipped with a linear order $\leq$ such that if $a \leq b$ then $a+c \leq b+c$, for every $a,b,c \in \Lam$. Then, a \emph{\Lam-metric space} is a pair $\ab{X,d}$ where $X$ is a set and $d \colon X \times X \to \Lam$ is a function satisfying the usual axioms for metric spaces. We can define intervals in \Lam\ and arcs and geodesic segments in $X$ analogously with the real case. A \emph{\Lam-tree} is then a \Lam-metric space $\ab{X,d}$ such that between any two points $x,y \in X$ there is a unique arc, denoted $[x,y]$, which is also a geodesic segment. When $\Lam=\R$ we obtain our familiar \R-trees, while setting $\Lam=\Z$ yields the class of graph-theoretic trees.

We can now ask analogues of the main question for ordered abelian groups other than \Lam. The notions of branchwise-real order tree and continuous \R-grading generalise readily to \emph{branchwise-\Lam\ tree order} and \emph{continuous \Lam-grading}. The main question then generalises to the following classification problem.

\begin{question}
	For which ordered abelian groups \Lam\ is it the case that every branchwise-\Lam\ tree order is continuously \Lam-gradable?
\end{question}

\Cref{res:non cont grad brto} provides a negative answer in the case $\Lam=\R$. On the other hand, an inductive argument shows that on every branchwise-\Z\ tree order we can indeed build a continuous \Z-grading. What about other ordered abelian groups?

A different line of questioning seeks to investigate the class of branchwise-real tree orders themselves. What properties can such tree orders possess? Let us take inspiration from the rigidity analysis of well-stratified trees (see \cite{10.2307/40378073} and references therein). Call a partial order \emph{rigid} if it has no non-trivial order-automorphisms. The following successively more difficult questions present themselves.

\begin{question}
	Does there exist a rigid branchwise-real tree order?
\end{question}

\begin{question}
	Does there exist a rigid branchwise-real tree order in which every branching node has the same degree?
\end{question}

\begin{question}
	Does there exist a rigid branchwise-real tree order in which every node is branching and has the same degree?
\end{question}

This line of investigation will be taken up in future work.

\section{Acknowledgements}

I wish to thank my supervisor, Joel David Hamkins, for much helpful discussion, guidance and proof-reading throughout the whole process of writing this paper. I would also like to thank the anonymous reviewer for helpful comments on a previous version. This research was supported by an EPSRC Studentship with project reference \emph{2271793}.

	\printbibliography

\end{document}